\documentclass[reqno]{amsart}

\usepackage{amssymb,amscd,color}

\theoremstyle{plain}
\newtheorem{theorem}{Theorem}[section]
\newtheorem{lemma}[theorem]{Lemma}
\newtheorem{proposition}[theorem]{Proposition}

\theoremstyle{definition}
\newtheorem{definition}[theorem]{Definition}
\newtheorem{problem}{Problem}

\theoremstyle{remark}
\newtheorem*{remark}{Remark}
\newtheorem*{remarks}{Remarks}

\newcommand{\N}{\hbox{\ensuremath{\mathbb{N}}}}
\newcommand{\Z}{\hbox{\ensuremath{\mathbb{Z}}}}

\newcommand{\R}{\hbox{\ensuremath{\mathbb{R}}}}
\newcommand{\C}{\hbox{\ensuremath{\mathbb{C}}}}
\newcommand{\E}{\mathcal {E}}
\newcommand{\LL}{\mathcal {L}}
\newcommand{\HH}{ \mathcal H}

\newcommand{\VV}{\mathcal V}
\newcommand{\e}{%
\mathfrak{E}}%
\newcommand{\Ggot}{%
\mathfrak{G}}%
\newcommand{\F}{\hbox{\ensuremath{\mathcal{F}}}}

\newcommand{\zo}{\ensuremath{\omega}}
\newcommand{\rank}{\text{rank}}
\newcommand{\sspan}{\text{span}}

\newcommand{\begeq}{\begin {equation}}

\newcommand{\SIS}{shift-invariant space} 
\newcommand{\SSIS}{Shift-invariant spaces}

\numberwithin{equation}{section}

\begin{document}


\title{Optimal Shift invariant spaces and their Parseval frame generators}


\author[Akram Aldroubi]{Akram~Aldroubi} \address{Department of Mathematics\\
Vanderbilt University\\ 1326 Stevenson Center\\ Nashville, TN 37240}
\email[Akram Aldroubi]{akram.aldroubi@vanderbilt.edu}

\author[Carlos Cabrelli]{Carlos~Cabrelli}
\address{Departamento de
Matem\'atica \\ Facultad de Ciencias Exactas y Naturales\\ Universidad
de Buenos Aires\\ Ciudad Universitaria, Pabell\'on I\\ 1428 Capital
Federal\\ ARGENTINA\\ and CONICET, Argentina}
\email[Carlos~Cabrelli]{cabrelli@dm.uba.ar} \thanks{The research of
Akram Aldroubi is supported in part by NSF Grant DMS-0504788. The research Douglas Hardin is supported in part
by NSF Grants DMS-0505756 and DMS-0532154. The research of
       Carlos Cabrelli and Ursula Molter is partially supported by
Grants: PICT 15033, CONICET, PIP 5650, UBACyT X058 and X108}

\author[Douglas Hardin]{Douglas~Hardin} \address{Department of Mathematics\\,
Vanderbilt University\\ 1326 Stevenson Center\\Nashville, TN 37240}
\email[Doug Hardin]{doug.hardin@vanderbilt.edu}

\author[Ursula Molter]{Ursula~Molter}
\address{Departamento de
Matem\'atica \\ Facultad de Ciencias Exactas y Naturales\\ Universidad
de Buenos Aires\\ Ciudad Universitaria, Pabell\'on I\\ 1428 Capital
Federal\\ ARGENTINA\\ and CONICET, Argentina}
\email[Ursula~M.~Molter]{umolter@dm.uba.ar}

\date{\today}

\begin{abstract}
Given a set of functions $\F=\{f_1,
\dots,f_m\} \subset L^2(\R^d)$, we study the problem of finding the shift-invariant space $V$ with $n$ generators $\{\varphi_1,\dots,\varphi_n\}$ that is ``closest'' to the functions of $\F$ in the sense that
  \begin {equation*}
V = \hbox{argmin}_{V' \in \mathcal V_n}\sum \limits_{i=1}^m w_i\|f_i-P_{V'}f_i\|^2,
\end {equation*}
where $w_i$s are positive weights, and $\mathcal V_n$ is the set of all shift-invariant spaces that can be generated  by $n$ or less generators. The Eckart-Young Theorem uses the singular value decomposition to provide a solution to a related  problem in finite dimension. We transform the problem under study into an uncountable set of finite dimensional problems each of which can be solved using an extension of the Eckart-Young Theorem. We prove that the finite dimensional solutions can be patched together and transformed to obtain the optimal shift-invariant space solution to the original problem, and we produce a Parseval frame for the optimal space.  
A typical application is the problem of finding a shift-invariant space model that describes a given class of signals or images (e.g., the class of chest X-Rays),   from the observation of a set of $m$ signals or images $f_1,
\dots,f_m$, which may be theoretical samples, or experimental data.  
\end{abstract}

\maketitle


\section{Introduction}
 In many signal and image processing applications, images and signals
are assumed to belong to some 
 \SIS\  of the form:
\begin{equation}
\label {sis}
\mathcal{S}(\Phi):=\hbox{closure}_{L_2} \;  \hbox{span} \{ \varphi_i(x-k): \;  i =1,\ldots, n, \, k \in \Z^d\} 
\end{equation}
       where $\Phi=\{\varphi_1,..., \varphi_n\}$ is a  set of functions in
$L^2(\R^d)$.        The  functions $\varphi_1 ,
\varphi_2, \ldots, \varphi_n$  are called a  {\em set of generators} for the
space $\mathcal S=\mathcal{S}(\Phi)=\mathcal{S}(\varphi_1,\ldots,\varphi_n)$  and any such space $\mathcal S$ is called a {\em finitely generated shift-invariant space (FSIS)} (see e.g.,
\cite {BDR94a, AG01}). For example, if
$n=1,\ d=1$ and
$\phi(x) = \hbox {sinc}(x)$, then the underlying space is the
space of \index{band-limited functions} band-limited functions (often
used in communications).

Finitely generated shift-invariant spaces,
 can have different sets of generators. The {\em length} of an FSIS  $\mathcal S$ is,
$$
\LL(\mathcal S)=\min  \{\ell \in \N: \exists \; \varphi_1,\ldots,\varphi_\ell \in \mathcal S \text{ with } \mathcal S=\mathcal S(\varphi_1,\ldots,\varphi_\ell)\}
$$
We will denote by $\mathcal V_n$  the set of all  shift-invariant invariant spaces with length less than or equal to $n$. That is, an element in $\mathcal V_n$ is a \SIS \  that has a set of $s$ generators with $s \leq n$.


In most applications, the shift-invariant space chosen to
describe the underlying class of signals is not
derived from experimental data - for example many signal
processing applications assume
``band-limitedness'' of the signal, which has theoretical advantages,
but generally does not necessarily reflect the underlying class of
signals accurately. 
Furthermore, in applications, the a priori
hypothesis that the class of signals
belongs to a shift-invariant space with a known number of generators, may
not be satisfied. For example,
the class of functions from which the data is drawn may not be a
shift-invariant space. Another example is when the shift-invariant space hypothesis is correct but
the assumptions about the number
of generators is wrong. A third example is when the a priori
hypothesis is correct but the data is
corrupted by noise. In addition, for computational considerations, a  shift-invariant space of length $m$ could be modeled by a shift-invariant model space with  length $n$ much smaller than $m$. 
For example, in Learning Theory, the problem of reducing the number of generators for a subspace of a reproducing kernel Hilbert space is also important for improving the efficiency and sparsity of
learning algorithms (see \cite{SZ04}).  In order to model classes of signals or images by FSIS in realistic cases, or to model a very large data set by a computationally manageable shift-invariant space, we consider the following problem:

\begin {problem} \label {P1}

Given a large set of  experimental data $\F = \{f_1, f_2, \ldots, f_m\}\subset L^2(\R^d)$, 
 we wish to determine a \SIS\   $V\in \mathcal{V}_n$ (where typically $n$ is chosen to be small compared to $m$) that models the signals in ``some'' best way.  
 For this purpose, we consider the following least squares problem:
  \begin {equation}
\label {opsp}
V = \hbox{argmin}_{V' \in \mathcal V_n}\sum \limits_{i=1}^m w_i\|f_i-P_{V'}f_i\|^2
\end {equation}
where $w_i$ are positive \index{weights} weights and where $P_{V'}$ is the
\index{orthogonal projection} orthogonal projection
on $V'$.

\end{problem}
A space $V$ satisfying \eqref{opsp} will be said to {\em solve Problem 1 for $(\F,w,n)$}.

The weights $w_i$ can
be chosen to normalize or to reflect our confidence about the data.
For example we can
choose $w_i=\|f_i\|^{-2}$ to place the data on a sphere or we can
choose a small weight $w_i$ for
a given  $f_i$ if ---due to \index{noise} noise or other factors---
our confidence
about the accuracy of $f_i$ is
low.  The goal is to see if we can perform operations on the observed data
 $\F = \{f_1, f_2,\ldots, f_m\}$ to construct (if it exists)
a shift-invariant space $\mathcal{S}(\Phi)$ whose length doesn't exceed a small number $n$, that minimizes the error with our data $\F$. 

Problem \ref {P1} can be
viewed as non-linear infinite dimensional constrained minimization
problem. It may also be viewed in light of the recent learning theory
  developed in \cite {BCDV05,CS02, SZ03}, and estimates of model fit in terms of
noise and approximation space may be derived. Beside the fundamental question of existence of an optimal space, it will be important for applications
to have a way to construct the generators of the optimal space if it exists, and to
estimate the error $\mathcal E(\F,w,n)= \sum \limits_{i=1}^m w_i\|f_i-P_{V}f_i\|^2$ where $ V\in \mathcal V_n$ is an optimal space for $\mathcal F, w$ and  $n.$

Typical  applications involve   large
 data sets (for
example consider the problem of finding a \SIS\ model for the collection of chest X-rays using data collected by a hospital  during the last 10 years).
The space
$\mathcal{S}(\F) $
generated by a set of experimental data contains all the data as
possible signals, but it is too large to be an appropriate model
for use in applications.  A space with a ``small'' number of
generators is more suitable, since if the space is chosen
correctly, it would reduce noise, and would give a computationally
manageable  model for a given application.

Least squares problems of the form above in {\em finite dimensional} spaces  can be solved using
the singular value decomposition (SVD).   \SSIS \, are infinite dimensional and the SVD cannot be applied directly.  However,  due to the special structure of \SIS s, the Fourier transform converts Problem \ref{P1} into finite dimensional least square problems at each frequency as will be discussed in Section \ref{sec-reduction}.



\section{Main Theorems}
\label {mainths}
In this paper we will sometimes deal with the standard Hilbert space $\C^N.$ Elements of this 
vector space are column vectors with $N$ coordinates.
 We will use the notation $A^t$ and $A^*$
 to denote the transpose and the conjugate transpose respectively of a complex matrix $A$.
 We will say that a vector $y \in \C^N$ is a left eigenvector of the matrix $A$ associated
 to the eigenvalue $\lambda$, if $y^t A = \lambda y^t$.

For clarity in the exposition, we will consider the unweighted case ($w_i=1, i=1,\ldots,m$).
The general case can be derived by simply applying the results of the unweighted case to the set of normalized observations $\F=\{f_1/w_1^2,\dots, f_m/w_m^2\}$.

The first theorem establishes the existence of an optimal space $V$. It also establishes  that $V$ can always be chosen to be a subspace of the shift-invariant space $\mathcal S (\F)$ generated by the totality of the data. This optimal space $V$ may not be  unique. However, under additional assumptions that are often satisfied in practice, there is only one optimal space $V$, as stated in Theorem \ref  {unique}.

\begin {theorem}
\label {existence}
Let $\F=\{f_1,\dots,f_m\}$ be a set of functions in $L^2(\R^d)$. Then
\begin{enumerate}
\item There exists $V\in \VV_n$ such that 
\begin{equation}\label{best}
\sum \limits_{i=1}^m \|f_i-P_Vf_i\|^2\le \sum \limits_{i=1}^m \|f_i-P_{V'}f_i\|^2,  \quad \forall \ V' \in \VV_n
\end{equation}
\item The optimal shift-invariant space  $V$ in \eqref{best} can be chosen  such that $V\subset \mathcal {S} (\F)$.
\end {enumerate}
\end {theorem}
\begin {remarks} 
\ 
\begin {itemize}
\item [(i)] Although we do not make the assumption that $n\le m$, if $n>m$, then $\mathcal S(\F)$ is an optimal space that belongs to $\mathcal V_n$. Thus, we will always assume that $n\le m$ for the remainder of this paper. 
\item [(ii)] In practice it will often be the case that $n$ is chosen (or found) to be much smaller than $m$. 
\end {itemize}
\end {remarks}
We still need to explicitly find an optimal space $V$ and estimate the error
\begin {equation}
\label {err}
\mathcal E (\F,n)=\min\limits_{V'\in \mathcal V_n}\sum \limits_{i=1}^m \|f_i-P_{V'}f_i\|^2.
\end {equation} 
To compute the error $\mathcal E (\F,n)$ we need to consider the Gramian matrix $G_{\F}$ of 
$\F=\{f_1,\dots,f_m\}.$ Specifically, the {\em Gramian} $G_\Phi$ of a set of functions $\Phi=\{\varphi_1,\dots,\varphi_n\}$ with elements in $L^2(\R^d)$ is defined to be the $n \times n$ matrix of $\Z^d$-periodic functions 

\begin{equation} \label{gram}
[G_{\Phi}(\zo)]_{i,j} = \sum_{k\in\Z^d}\widehat \varphi_i(\zo+k)\overline{\widehat \varphi_j(\zo+k)}, \qquad \zo\in \R^d,
\end{equation}
where $\widehat \varphi_i$ denotes the Fourier transform of $\varphi_i$, and where $\overline{\widehat \varphi_i}$ denotes the complex conjugate  of $\widehat \varphi_i$. It is known that $G_\Phi$ is $\Z^d$-periodic non-negative and self-adjoint for almost every $\zo$. In this paper, we use the following definition for the Fourier transform of a function $\phi \in L^2(\R^d)$:
\begin{equation}
\hat \phi (\zo):=\int_{  \R^d} \phi(x) e^{-i2\pi x^t \zo}\, dx, \qquad \zo\in \R^d, 
\end{equation}
where $dx$ denotes Lebesgue measure on $\R^d$. 

If $V = \{v_1, \dots, v_n\}$ is a set of vectors in a Hilbert space $\HH$, we will denote by
$\Ggot(V) = \Ggot(v_1, \dots, v_n)$ the matrix 
\begin{equation} \label{gram-1}
[\Ggot(V)]_{i,j} =  \langle v_i, v_j \rangle_{\HH} \quad i, j = 1, \dots, n.
\end{equation}

Our next  theorem produces a generator for an optimal space $V$ and   provides a formula for the exact value of  the error, but we first recall the definition and some properties of  frames used in its statement (see for example \cite {CC97,HL00,HLW02}). 
 \begin{definition} 
 \label {framedef}
 Let $\HH$ be a Hilbert space and $\{u_i\}_{i\in I}$ a countable subset of $\HH$. The set $\{u_i\}_{i\in I}$ is said to form a {\em frame}  for $\HH$ if there exist $q,Q>0$ such that 
 \[
 q\|f\|^2\le \sum\limits_{i\in I}|<f,u_i>|^2\le Q \|f\|^2, \quad \forall \, f \in \HH.
 \]
 If $q=Q$, then $\{u_i\}_{i\in I}$ is called a {\em tight frame}, and it is called a {\em Parseval frame} if $q=Q=1$.
 \end {definition}
 If $\{u_i\}_{i\in I}$ is a Parseval frame for a subspace $W$ of a Hilbert space $\HH$, and if $a \in \HH$, then the orthogonal projection of  $a$ onto $W$ is given by:
  \begin{equation} \label{proj}
P_W(a) = \sum_{i\in I} \langle a, u_i \rangle u_i.
\end{equation}
Thus, a Parseval frames acts as if it were an orthonormal basis of $W$, even though it may not be one.
\begin {theorem} \label{main} Under the same assumptions as in Theorem \ref {existence}, let   $\lambda_1(\zo)\ge\lambda_2(\zo)\ge\dots\ge \lambda_m(\zo)$ be the eigenvalues of the Gramian $G_{\F}(\zo)$.  Then
\label {construction}
\begin {enumerate}
\item The eigenvalues $\lambda_i(\zo) $, $1\le i \le m$ are $\Z^d$-periodic, measurable functions in $L^2([0,1]^d)$ and 
\begin{equation} 
\label{err}
\mathcal E(\F,n)= \sum \limits_{i=n+1}^m\int\limits_{[0,1]^d}\lambda_i(\zo)d\zo .
\end{equation}

\item \label{teo-4} Let $E_i:=\{\zo\;:\; \lambda_i(\zo)\ne 0\}$, and define $\tilde \sigma_i(\zo)=\lambda_i^{-1/2}(\zo)$ on $E_i$ and $\tilde \sigma_i(\zo)=0 $ on $E_i^c$. Then, there exists a choice of measurable left eigenvectors   $y_1(\zo),\dots,y_n(\zo)$ with $y_i=(y_{i1},...,y_{im})^t, i=1,...,n,$ associated with the first $n$ largest eigenvalues of $G_{\F}(\zo)$ such that the functions defined by
\[
\hat \varphi_i(\zo)=\tilde \sigma_i(\zo)\sum \limits_{j=1}^my_{ij}(\zo)\hat f_j(\zo), \quad i=1,\dots,n, \; \zo\in \R^d
\] 
are  in $L^2(\R^d)$. Furthermore, the corresponding set of  functions $\Phi=\{\varphi_1,\dots,\varphi_n\}$ is a generator for an optimal space $V$ and the set $\{\varphi_i(\cdot-k), k \in \Z^d, i=1,\dots,n\}$ is a Parseval frame for $V$. 
\end{enumerate}
 \end {theorem}
 
 The following example shows that the optimal space $V$ does not need to be unique.
Let $m=2$, $n=1$, and  let $f_1,f_2$ be two orthonormal functions. For this situation, $G_{\F}(\zo)$ is the $2\times 2$ identity matrix for almost all $\zo \in \R^d$.   It follows that any function $\varphi=c_1f_1+c_2f_2$ with $c=(c_1,c_2)$ a unit vector  in $\R^2$ generates an optimal space and $\mathcal E(\F,1)=1$. Obviously, in this particular case there are infinitely many optimal spaces.  However, under some mild assumptions, there exists a unique optimal space $V$ as described in the following theorem:
 \begin {theorem}
 \label {unique}
 Let $\F=\{f_1,\dots,f_m\}$  functions  in $L^2(\R^d)$ be given.  If $\lambda_n(\zo)> \lambda_{n+1}(\zo)$ for almost all $\zo$, then the optimal space $V$ in \eqref{best} is unique. In this case,  $n\le r_{\min}=\min\limits_{\zo\in [0,1]^d} \mathrm{rank}\,G_{\F}(\zo)$ and  the set  $\{\varphi(\cdot-k), k \in \Z^d, i=1,\dots,n\}$ in part (\ref{teo-4}) of Theorem \ref {construction} is an orthonormal basis for $V.$ 
\end {theorem}
\begin {remarks}
\ 
\begin {itemize}
\item[(i)] In case that $n = \mathcal L(\mathcal {S}(f_1,...,f_m)),$ Theorem \ref{main} gives a proof of the
known result that every FSIS has a set of generators forming a Parseval frame.
\item[(ii)] It will be clear from the proofs of Theorems \ref{existence} and \ref{main} that
the optimal space $V$ can be decomposed as $V=S(\varphi_1)\bigoplus...\bigoplus S(\varphi_{\ell})$
where $\ell = \mathcal L(V),$ the direct sum is orthogonal and each $\varphi_i$ is a Parseval frame generator of $S(\varphi_i).$
\item [(iii)]Theorem \ref {unique} can only be used when $n <m$. When $n=m$ then $\mathcal S (\F)$ is an optimal space and  it is the unique optimal space if and only if $\LL(\mathcal S (\F))=m$. 
\item [(iv)] Obviously, if $n=m$ then the error between the model and the
observation is null. However, by plotting the error in \eqref{err} in terms of the number of generators, an optimal
number $n$ may be heuristically derived.  Alternatively, one may choose   $n$ so that a cost functional (depending on the error and on $n$) is optimized as in other
dimension reduction schemes.
\end{itemize}
\end {remarks}

\section{Preliminaries on Finitely Generated Shift-Invariant Spaces}

\label{section-sis}

In this section we state some known results about finitely generated shift-invariant spaces that we will need later.
See for example (\cite{Hel64, BDR94, BDR94a, RS95, Bow00}.)

 We need first to introduce some definitions.

Given $f \in L^2(\R^d)$ and $x \in \R^d$ the {\em fiber} of $f$ at $x$  
is the sequence $\Gamma_x f = \{f(x+k): k \in \Z^d\}.$

If $V$ is a FSIS (recall Definition \eqref {sis}) and $\omega \in [0,1]^d$ we set  
$V_\omega = \{\Gamma_\omega \hat f ;  f \in V\}$ the {\em fiber space} associated to $V$ and $\omega$.

If $\mathcal M$ is a closed subspace of a Hilbert space $\mathcal H$, throughout this article  we will denote by $P_{\mathcal M}$  the orthogonal projection operator  in $\mathcal H$ onto $\mathcal M.$ 

With this notation we have:

\begin {lemma} \label {OpGam} 
If $f\in L^2(\R^d)$, then 
\begin{enumerate}
\item  The sequence $(\Gamma_\zo\hat f)_k=(\hat f(\zo+k))$ is a well-defined sequence in $\ell_2(\Z^d)$ a.e. $\zo \in \R^d$; and
\item $\|\Gamma_\zo\hat f\|_{\ell_2}$ is a measurable function of $\zo$ and $\|f\|^2=\|\hat f\|^2=\int\limits_{[0,1]^d} \|\Gamma_\zo\hat f\|^2_{\ell_2}d\zo.$
\end{enumerate}
\end {lemma}

\begin{lemma}\label{fibers}
Let $V$ be a FSIS in $ L^2(\R^d).$ Then we have:
\renewcommand{\labelenumi}{{\normalfont \roman{enumi}.}}
\renewcommand{\labelenumii}{{\normalfont \alph{enumii}.}}
\begin{enumerate}
\item
$V_{\omega}$ is a closed subspace of $\ell_2(\Z^d)$ for almost all $\omega \in [0,1]^d.$
\item
$V=\{ f \in  L^2(\R^d) :  \Gamma_\omega \hat f \in V_{\omega}  \text{ for almost all } \omega \in [0,1]^d\}.$
\item
 For  each $f \in L^2(\R^d)$ we have that $||\Gamma_{\omega} (\widehat{P_Vf})||_{\ell_2}$ is a measurable function of the variable $\omega$ and
 $\Gamma_{\omega} (\widehat{P_Vf}) = \Gamma_{\omega}P_{\hat V} \hat f = P_{V_{\omega}}(\Gamma_{\omega}\hat f).$
\item
Let  $\varphi_1,...,\varphi_r \in L^2(\R^d)$. We have  that 
\begin{enumerate} 
\item
$\{\varphi_1,...,\varphi_r\}$ is a set of generators of $V,$  if and only if 
 the fibers \\ $\Gamma_{\omega} \hat \varphi_1,...,\Gamma_{\omega} \hat \varphi_r$ span $V_{\omega}$ for almost all $\omega \in [0,1]^d$
 \item
 the integer translates of $\varphi_1,...,\varphi_r$ are a frame  of $V,$  if and only if \\  $\Gamma_{\omega} \hat \varphi_1,...,\Gamma_{\omega} \hat \varphi_r$ are a 
frame of $V_{\omega}$ with the same frame bounds, for almost all  $\omega \in [0,1]^d.$
\end{enumerate}
\end{enumerate}
\end{lemma}

\begin {lemma}
\label {GramPer}
Let $\F=\{f_1,\dots,f_m\}$ be functions  in $L^2(\R^d)$ and let $A(\zo)$ be the infinite matrix $A_{kj}(\zo)=\left(\Gamma_\zo \hat f_j\right)(k)=\hat f_j(\zo+k)$, $j=1,\dots,m$, $k\in \Z^d$, and $\zo\in \R^d$. Then $G_{\F}(\zo)=A^t(\zo)\overline {A(\zo)}$, and $\rank\,G_{\F}(\zo)=\rank\, A(\zo)=\rank\, A^*(\zo), \;  a.e.\;  \zo \in \R^d$. In particular, $G_{\F}(\zo) = \Ggot(\Gamma_\zo \hat f_1, \dots, \Gamma_\zo \hat f_m)$.
\end {lemma}


\section{Proofs}

To prove the theorems in Section \ref {mainths}, we proceed in several steps. First we reduce the optimization problem into an uncountable set of finite dimensional problems in the Hilbert space $\mathcal H= \ell_2(\Z^d)$. We then apply the Eckart-Young Theorem to prove that the reduced problems have solutions. Finally, we construct the generators of the optimal space patching together the  solutions of the reduced problems to obtain the solution to the original problem.

\subsection{Reduction} \label{sec-reduction}
In this section, we reduce Problem \ref {P1} to a set of finite dimensional problems.
To see this let us first consider the following : 


\begin{problem} \label{P2}

Let $H$ be a Hilbert space, $n$, $m$ positive integers and 
$A= \{a_1,...,a_m \}$ a set of vectors in $ \HH$. We want to find a closed subspace $S$ of $\HH$ with $\text{dim}(S)\leq n$ 
that satisfies
\begin{equation}\label{min-I}
\sum\limits_{i=1}^m  \| a_i-P_S a_i  \|^2\leq \sum\limits_{i=1}^m\|a_i-P_{S'}a_i \|^2
\end{equation}
for every subspace $S'\subset \HH $ with $\text{dim}(S')\leq n.$

If such an $S$ exists, we say that $S$ 
{\em solves Problem \ref{P2} for the data $(A,n)$}.

If $B=\{b_1,...,b_r\}$ is a set of vectors from $\HH$ with $S= \text{span}(B)$ we will say that
the vectors in {\em $B$ solve Problem \ref{P2} for the data $(A,n)$}.
The error for Problem \ref{P2} is  
$$\e(A,n) =\min\limits_{\text{dim}(S')\leq n} \sum\limits_{i=1}^m\|a_i-P_{S'}a_i \|^2.$$
\end{problem}
Note that in Problem \ref{P2} we take the minimum over all subspaces of dimension less than $n$, while
in Problem \ref{P1} the minimization is taken over a particular class of infinite dimensional subspaces,
so the two problems are essentially different.

In the next section we  state and prove an extension of the Eckart-Young theorem. We conclude from this extension that  Problem \ref{P2}  always has a solution
for any set of data $(A,n)$ in an arbitrary Hilbert space. That is, given $A$ and $n$ there always exists a subspace $S$ with dim($S$) $\leq n $ satisfying (\ref{min-I}).
We will also see  that a solution $S$ can be chosen in such a way that $S \subset \text{span}(A)$ when $n \leq \text{dim(span}(A))$.

Before proving these results let us see how Problem \ref{P2} helps our original question.

Let  $\mathcal{F}=\{f_1,\ldots, f_m\}\subset L^2(\R^d)$.  We want to find out if there exists $V\in \mathcal{V}_n$ such that $V$ minimizes $\sum_{i=1}^{m} \|f_i - P_{V}f_i\|^2$.  Using Lemma \ref {OpGam} we obtain that for any $V \in \mathcal{V}_n,$
\begin{align}
\sum_{i=1}^m \|f_i-P_Vf_i\|^2&=\sum_{i=1}^m \int_{[0,1]^d} \|\Gamma_\omega\hat f_i -\Gamma_\omega\widehat{P_Vf_i}\|_{\ell_2}^2\, d\omega \notag 
\\
&=\int_{[0,1]^d} \sum_{i=1}^m \|\Gamma_\omega\hat f_i -\Gamma_\omega\widehat{P_Vf_i}\|_{\ell_2}^2\, d\omega.
 \label{mintransform}
 \end{align}

By Lemma \ref{fibers}(iii), $\Gamma_\omega\widehat{P_Vf_i} = P_{V_\omega}\Gamma_\omega\hat f_i.$
So from (\ref{mintransform}) we conclude that,
\begin{equation}\label{min-II}
\sum_{i=1}^m \|f_i-P_Vf_i\|^2 = \int_{[0,1]^d} \sum_{i=1}^m \|\Gamma_\omega\hat f_i -P_{V_\omega}\Gamma_\omega \hat f_i \|_{\ell_2}^2\, d\omega.
\end{equation}

The sum inside the integral on the right hand-side of (\ref{min-II})  is of the same type than the sum 
that is involved in Problem \ref{P2} in the case that $\HH = \ell_2(\Z^d)$ and $S=V_{\omega}$.
Since we are assuming that Problem \ref{P2} always has a solution, we know that for almost each $\omega \in [0,1]^d$ there exists a subspace $S_\omega\subset \ell_2(\Z^d)$ that solves Problem \ref{P2} for the data 
$(\F_{\omega},n)$ where 
$\F_{\omega}=\{\Gamma_\omega\hat f_1,...,\Gamma_\omega\hat f_m\}.$ Note that the subspace
 $S_\omega$ does not need to be related with the fiber space of any FSIS.
If the function 
$\omega \longmapsto  \sum_{i=1}^m \|\Gamma_\omega\hat f_i -P_{S_\omega}\Gamma_\omega \hat f_i \|_{\ell_2}^2\,$ were a measurable function of $\omega$ then we would have

\begin{equation}
\int_{[0,1]^d} \sum_{i=1}^m \|\Gamma_\omega\hat f_i -P_{S_\omega}\Gamma_\omega \hat f_i \|_{\ell_2}^2\, d\omega \leq \sum_{i=1}^m \|f_i-P_{V'}f_i\|^2
\end{equation}
for every $V' \in \VV_n.$

Therefore, in case that there exists  a FSIS   $V \in \VV_n$ such that
 $V_\omega = S_\omega$ a.e. $\omega \in [0,1]^d$, then by Lemmas \ref{OpGam} and \ref{fibers} the above function would be measurable and $V$ necessarily will be a solution  to Problem \ref{P1}, since
 \begin{multline} \label{min}
 \sum_{i=1}^m \|f_i-P_Vf_i\|^2 =  \int_{[0,1]^d} \sum_{i=1}^m \|\Gamma_\omega\hat f_i -P_{S_\omega}\Gamma_\omega \hat f_i \|_{\ell_2}^2\, d\omega  \leq \\
 \int_{[0,1]^d} \sum_{i=1}^m \|\Gamma_\omega\hat f_i -P_{V'_\omega}\Gamma_\omega \hat f_i \|_{\ell_2}^2\, d\omega  =\sum_{i=1}^m \|f_i-P_{V'}f_i\|^2 
\end{multline}
for every $V' \in \VV_n.$

 We will see later that such a FSIS indeed exists. More precisely we will construct a set
 of generators such that its integer translates form a frame of the optimal FSIS. We will do that by patching together the fibers of the generators of each of the optimal subspaces $S_\omega.$

                                 
\subsection{Solution to Problem \ref{P2}}

We now prove that Problem \ref{P2} always has a solution. 

\begin {theorem}\label{Eckart-H}
Let $\HH$ be an infinite dimensional  Hilbert space, $\F=\{f_1,...,f_m\} \subset \HH$, 
$\mathcal X = \text{span }\{f_1,...,f_m\},$  
 $\lambda_1 \geq ...\geq \lambda_m$ the eigenvalues of the matrix $\Ggot(\F)$ defined as in \eqref{gram-1}  and $y_1,...,y_m \in \C^m$, with
 $y_i=(y_{i1},...,y_{im})^t$ orthonormal  left eigenvectors associated to the eigenvalues
 $\lambda_1,...,\lambda_m.$ Let $r = \dim {\mathcal X} =\text{\rm rank}\ \Ggot(\F) $.

Define the vectors $q_1,...,q_n \in \HH$ by
 \begin {equation}
 \label {Eckart-H1}
 q _i =
 \tilde \sigma_i \sum_{j=1}^m y_{ij}f_j,   \;\;\;\;\; \quad i=1,...,n \\
  \end{equation}
  where $\tilde \sigma_i = \lambda_i^{-1/2}$ if $\lambda_i \neq 0$, and $\tilde \sigma_i = 
  0$ otherwise.
Then  $\{q_1,...,q_n\}$ is a Parseval frame of $W= \text{span }   \{q_1,...,q_n\}$
and the subspace $W$ is optimal in the sense that

$$
\e(\F,n)=\sum_{i=1}^m \|f_i -P_{W}f_i\|^2
 \le \sum_{i=1}^m \|f_i-P_{W'}f_i\|^2, \quad \forall \text{ subspace } \; W',  \dim W'\le n.
$$

Furthermore we have the following formula for the error
\begin{equation}\label{error1}
\e (\F,n)=  \sum_{i=n+1}^m \lambda_i.
\end{equation}
\end{theorem}

\begin{remark} If $r$ is small (i.e. $r \leq n$) then all the vectors $q_{r+1}, \dots, q_n$ are null and $\{q_1,...,q_r\}$ is
an orthonormal set. 

One could also choose $q_{r+1}, \dots, q_n$ to be any orthonormal set in the orthogonal complement of $\mathcal X$ and so obtain an orthonormal set of $n$ elements and the formula for the error would still hold.
\end{remark}

If $\HH$ is finite dimensional and $n \leq r$, then Theorem \ref{Eckart-H} is a consequence of the Eckart-Young theorem (see appendix).
To prove Theorem \ref{Eckart-H} we will  reduce it to the finite dimensional case and then  use the Eckart-Young result.

We first need the following Lemma:


\begin{lemma}
\label {reduction}
Let $\HH$ be a Hilbert space, $\{f_1,\cdots,f_m\}\subset \HH$, $\mathcal X= \sspan \{f_1,\cdots,f_m\}$. Assume that there exists $M \subset \HH$ with $\dim M\le n$ such that
\[
\sum\limits_{i=1}^m\|f_i-P_Mf_i\|^2\le \sum\limits_{i=1}^m\|f_i-P_{M'}f_i\|^2
\]
for any subspace $M' \subset \HH$ with $\dim M'\le n$, then there exists $W \subset \mathcal X$, with $\dim W\le n$, such that 
\[
\sum\limits_{i=1}^m\|f_i-P_Wf_i\|^2=\sum\limits_{i=1}^m\|f_i-P_Mf_i\|^2.
\]
\end{lemma}

\begin{proof}
 Define the subspace $W=P_\mathcal XM$ as the orthogonal projection of $M$ onto 
$\mathcal X$. By construction, $W\subset \mathcal X$, and $\dim W \le n$.

 Let $f \in \mathcal X$, then we have
\begin{eqnarray*}
\|f-P_Wf\|^2 &=& \inf \{\|f-g\|^2 : \; g \in W\} \\
& \le & \|f-P_\mathcal XP_Mf\|^2 \\
& = & \|P_\mathcal Xf-P_\mathcal XP_Mf\|^2 \\
& \le & \|f-P_Mf\|^2.
\end{eqnarray*}
\end{proof}
This Lemma shows that, in a possibly infinite dimensional Hilbert space $\HH$, the problem of finding a finite dimensional subspace  $M \subset \HH$ with $\dim M\le n$ that ``best approximates'' $m$ vectors $\{f_1,\dots,f_m\}$,
can always be reduced to a search in the finite dimensional space $\mathcal {X}=\sspan \{f_1,\cdots,f_m\}$. 

We  now prove Theorem \ref{Eckart-H}.
\begin{proof} (of Theorem \ref{Eckart-H}).
Let $\tau: \mathcal X \longrightarrow \C^m$ be an isometric isomorphism.
Set $b_i=\tau(f_i), $ and let $B$ be the matrix having the vectors $b_i$ as columns.
 So,  $r = \dim \mathcal X = \text{rank}(B)$ and 
$B^t\overline B$ coincides with $\Ggot(\F) = \{<f_i,f_j>_{\HH}\}_{i,j}$.

Choose  orthonormal left eigenvectors $y_1,...,y_m \in \C^m$, with
 $y_i=(y_{i1},...,y_{im})^t$   associated to the eigenvalues
 $\lambda_1\geq...\geq\lambda_m$ of $B^t\overline B$, and  define the vectors
 \begin {equation}
  u _i =
 \tilde \sigma_i \sum_{j=1}^m y_{ij}b_j,   \quad \quad  i=1,...,n \\
  \end{equation}
  where as before $\tilde \sigma_i = \lambda_i^{-1/2}$ if $\lambda_i \neq 0$, and $\tilde \sigma_i = 
  0$ otherwise. 
  
  Then, if $n \leq r$ by Theorem \ref{E-Y} in appendix,
  the subspace $M \subset \C^m, M= \text{span } \{u_1,...,u_n\}$ satisfies:
  
\begin{equation}\label{xx}
\sum\limits_{i=1}^m\|b_i-P_Mb_i\|^2\le \sum\limits_{i=1}^m\|b_i-P_{M'}b_i\|^2
\end{equation}
  for every subspace $M' \subset \C^m$
 with dim $M'\leq n.$ 
 
 If however, $n \geq r$ then the left side of \eqref{xx} is $0$ and therefore the inequality is also satisfied. 

Setting  $W=\tau^{-1}(M)$ and noting that $\tau^{-1}P_M=P_W$ we have from (\ref{xx})
 \begin{equation}
\e(B,n)=\e(\F,n)=\sum\limits_{i=1}^m\|f_i-P_Wf_i\|^2\le \sum\limits_{i=1}^m\|f_i-P_{W'}f_i\|^2
\end{equation}
for every subspace $W'\subset \HH$ with dim $W' \leq n.$
So, $W\subset \HH$ is optimal for $(\F,n)$ and $q_i=\tau^{-1}(u_i),\:\;\; i=1,...,n$ is a Parseval frame for $W.$
Furthermore, the  formula (\ref{error1}) also holds.
\end{proof}
\begin{remarks}
\ 

\label {RemSVD}
\begin {itemize}
\item [(i)] If $n>m$ the optimal space $W$ is not unique since any space $W'$ of dimension $n$ containing $\sspan\{a_1, \dots, a_m\}$  will be optimal. The same argument also shows that  the space $W$ is not unique if $n> r= \text{dim } \mathcal X$.
\item [(ii)]If $n \leq r$, the vectors $u_i$ and $y_i$ are related by $\sqrt {\lambda_i} u_i= Ay_i$ as described in the Appendix.
\end {itemize}
\end {remarks}

\subsection{Solution to Problem \ref{P1}}

In order to solve Problem \ref{P1}, we need the following technical proposition 
concerning the measurability of the eigenvalues and the existence of measurable
eigenvectors of a non-negative matrix with measurable entries (cf. \cite[Lemma 2.3.5]{RS95}.)
 \begin {lemma}
\label {diagonalization}
Let $G=G(\zo)$ be an $m\times m$ self-adjoint matrix of measurable functions defined on  a   measurable subset  $E\subset \R^d$ with eigenvalues $\lambda_1(\zo)\ge\lambda_2(\zo)\ge\dots\ge\lambda_m(\zo)$.  Then the eigenvalues $\lambda_i$, $i=1,\ldots, m$, are measurable on $E$ and there exists an $m\times m$  matrix of measurable functions $U=U(\zo)$   on $E$  such that $U(\zo)U^*(\zo)=I$ a.e. $\zo \in E$ and such that
\begin{equation} \label{Gfac}
G(\zo)=U(\zo)\Lambda(\zo)U^*(\zo) \; \quad a.e. \; \zo \in E 
\end{equation}
where   $\Lambda(\zo):=diag(\lambda_1(\zo),\dots,\lambda_m(\zo))$.
\end {lemma}

\begin{proof} (of Theorems \ref{existence} and \ref{main})

In what follows we will apply Theorem \ref{Eckart-H} to find the solution to Problem \ref{P1}.
As before, let $\F=\{f_1,\dots,f_m\} \subset L^2(\R^d)$ and for $\omega \in [0,1]^d$ let  $G_{\F}(\omega)$ be the associated Gramian matrix with eigenvalues $\lambda_1(\omega) \geq ...\geq \lambda_m(\omega)\ge 0$.  Let $U(\zo)$   be a measurable  $m\times m$ matrix as in Lemma~\ref{diagonalization}. Since $G_{\F}(\omega)$ is $\Z^d$-periodic on $\R^d$, we can  choose
$U(\zo)$ to be $\Z^d$-periodic as well.   Let $U_i(\zo)$ denote  the $i$-th  row of $U(\zo)$.  Multiplying (\ref{Gfac})
on the left by  $U^*(\zo)$  shows that   $y_i(\zo):=
U_i(\zo)^*$ is a left-eigenvector of $G(\zo)$ with eigenvalue $\lambda_i(\zo)$ for $i=1,\ldots, m$.  Furthermore,  
the left eigenvectors  $y_i(\zo)=(y_{i1}(\zo),...,y_{im}(\zo))^t$, $i=1, \ldots, m$,  form an orthonormal basis of $\C^m$.

For each fixed $\zo\in [0,1]^d$, we consider Problem \ref{P2} in the space $\ell_2(\Z^d)$ for the data $(\F_\omega,n)$ with 
$\F_{\omega}=\{\Gamma_\omega\hat f_1,...,\Gamma_\omega\hat f_m\}$.
Define  $q_1(\zo),...,q_n(\zo) \in \ell_2(\Z^d)$  by
 \begin {equation}	 \label {Eckart-l2}
 q _i(\zo) = \tilde \sigma_i(\zo) \sum_{j=1}^m y_{ij}(\zo) \Gamma_{\zo} \hat f_j,   \;\;\;\;\; \quad i=1,...,n   
 \end{equation}
where   $\tilde \sigma_i(\zo) = \lambda_i^{-1/2}(\zo)$ if $\lambda_i(\zo) \neq 0$, and $\tilde \sigma_i(\zo) =   0$ otherwise. Since
$G_{\F}(\omega) = \Ggot(\F_{\omega})$ (see  Lemma~\ref{GramPer}),
 Theorem \ref{Eckart-H} shows that
 the space $S_{\omega}:= \text{span }   \{q_1(\zo),...,q_n(\zo)\}$ optimizes Problem \ref{P2}.  Moreover, the  vectors  $\{q_1(\zo),...,q_n(\zo)\}$ form a Parseval frame for   
$S_{\omega} $ and we have the following formula for the error
\begin{equation} \label{error}
\e(\F_{\zo},n)=  \sum_{i=n+1}^m \lambda_i(\zo).
\end{equation}

Define now the functions $h_i:\R^d\rightarrow \C$, $i=1,...,m$, 
\begin{equation}\label{measurable}
h_i(\zo) =  \tilde \sigma_i(\zo) \sum_{j=1}^m y_{ij}(\zo) \hat f_j(\omega).
\end{equation}
 
Since $\tilde \sigma_i$ and $y_i$ are measurable  functions of $\zo$, then  $h_i$ is also   measurable.
 Moreover, $h_i$  is in $L^2(\R^d)$ as the following simple argument shows.
 Since
 \begin{equation*}
 |h_i(\zo)|^2 = h_i(\zo) \overline{h_i}(\zo) =
 \tilde \sigma_i(\zo)^2 \sum_{j,s=1}^m y_{ij}(\zo) \hat f_j(\zo) \overline {\hat f_s}(\zo) \;
  \overline {y}_{is}(\zo)
 \end{equation*}
 we have (using that if $y_i$ is a left eigenvector of the self-adjoint matrix $G_{\F}$, then $\overline {y_i}$ is a right eigenvector for that matrix associated to the same eigenvalue),
 \begin{align} \label{xxx}
 \sum_{k \in \Z^d} |h_i(\zo+k)|^2 = & 
 \tilde \sigma_i(\zo)^2  \sum_{j=1}^m y_{ij}(\zo) \sum_{s=1}^m [G_{\F}(\zo)]_{js}   \overline {y}_{is}(\zo) \notag \\
=&\tilde \sigma_i(\zo)^2 \lambda_i(\zo) \sum_{j=1}^m y_{ij}(\zo) \overline {y}_{ij}(\zo) =
   \tilde \sigma_i(\zo)^2 \lambda_i(\zo).
 \end{align}
 If $\lambda_i(\zo) \neq 0$ then the product in (\ref{xxx}) is one, otherwise it is zero.
 That is  
 $\sum_{k \in \Z^d} |h_i(\zo+k)|^2= \mathbf{1}_{\{\zo:\lambda_i(\zo) > 0 \}}$ and by
  Lemma \ref{OpGam},  $\|h_i\| \leq 1.$
 
 Now define functions $\varphi_1,\dots,\varphi_n$ in $L^2(\R^d)$ by:
 $$\hat \varphi_i(\zo) = h_i(\zo), \quad i=1,\dots,n$$
and let $V=S(\varphi_1,\dots,\varphi_n).$   The space $V$ is a shift invariant space of
length no bigger than $n$. So $V \in \mathcal V_n.$
Furthermore, by Lemma \ref{fibers} (iv-a), the space $V_{\zo}$ is spanned by $\Gamma_{\zo} \hat \varphi_i, i=1,...,n.$

Since $(\Gamma_{\zo}\hat \varphi_i)(k) = h_i(\zo+k) =q_i(\zo)(k), k \in \Z^d,  i=1,\dots,n$ a.e.,
then $V_{\zo} = S_{\zo}$ (the optimal space for the data ($\F_{\zo},n))$ in $\ell_2(\Z^d)$.

By equation (\ref{min}) and the comment before, $V$ is optimal, that is $V$ solves Problem \ref{P1} for the data $(\F,n)$.

 Now, since $\{\Gamma_{\zo} \hat \varphi_1,...,\Gamma_{\zo} \hat \varphi_n\}$ is a Parseval frame of 
 $S_{\zo}$ for a.e. $\zo \in [0,1]^d$ then  by Lemma \ref{fibers} (iv-b) the integer translates
  of $\varphi_1,...,\varphi_n$ form a Parseval frame of $V.$
On the other hand, Formula (\ref{min-II}) says that 
\begin{equation}
\mathcal E (\F,n) = \int_{[0,1]^d} \e(\F_{\zo},n)\; d\zo. 
\end{equation}
Thus using (\ref{error}) we have that $\E(\F,n) = \sum_{i=n+1}^m  \int_{[0,1]^d} \lambda_i(\zo) \; d\zo.$
 \end{proof}
 
\begin{proof} (of Theorem \ref{unique}).

Under the hypothesis of Theorem \ref{unique}, Theorem \ref{Eckart-H} and the Remark after it guarantee the uniqueness of the optimal spaces $S_\zo$ associated to the data
($\F_{\omega},n$)
 for almost all $\zo$. Since these fiber spaces characterize the optimal space $V,$ 
then the theorem follows.
\end{proof}

 \section*{Appendix}
 
 \subsection{Best linear approximation and the SVD}

\label{sec-svd}

Here we review the singular value decomposition (SVD) of a matrix and  
its relation to  finite dimensional  least-squares problems.  For an overview see \cite {Ste93}, and for a very detailed treatment see for example \cite{HJ85}.

We  start with the following proposition.
 \begin{proposition}[SVD]
 Let $A = [a_1, a_2, \dots, a_m]$ be the matrix with columns $a_i \in \C^N$,  $m \leq N$.
Let  $r:=\dim \text {span } \{
 a_1,\ldots, a_m \}$. 
 Then  there are $m$ numbers $\lambda_1\ge \lambda_2\ge \cdots \lambda_r>\lambda_{r+1}=\cdots =\lambda_m=  0$,  an  orthonormal 
 collection of $m$ (column) vectors $y_1,\ldots, y_m\in \C^{m}$, and
 an orthonormal  collection of $m$ (column) vectors $u_1, \ldots, u_m\in \C^N$
 such that 
 \begin{equation}
 \label{SVD}
 A=\sum_{k=1}^m \sqrt{\lambda_k} u_k y_k^* = U\Lambda^{1/2} Y^* 
  \end{equation}
  Where $U \in \C^{N\times m}$ is the matrix $U = [u_1, \dots, u_m],
 \;\Lambda^{1/2} = \text{diag} (\lambda_1^{1/2},...,\lambda_m^{1/2}),$
 and  $Y = \{y_1,...,y_m\} \in C^{m \times m}$
  with $U^*U=I_m = Y^*Y=YY^*.$

  The representation of $A$ given in (\ref{SVD}) is called {\em the singular value decomposition} (SVD) of $A$. 
 \end{proposition}

 The SVD of a matrix  $A$ can be obtained as follows.  Consider the matrix
  $A^*A \in \C^{m\times m}$. Since $A^*A$ is self-adjoint and positive semi-definite, 
 its eigenvalues $\lambda_1\ge \lambda_2\ge \cdots \ge \lambda_m$   are nonnegative and the associated eigenvectors $y_1,\ldots, y_m$ can be chosen to form an orthonormal basis of $\C^m$.   Note that the rank $r$ of $A$ corresponds to the largest index $i$ such that $\lambda_i>0$. 
  The left singular vectors $u_1,\ldots, u_r$ can then be obtained from 
$$\sqrt{\lambda_i}u_i=Ay_i, \;\;  \text{ that is } u_i = \lambda_i^{-1/2} \sum _{j = 1}^m y_{ij} a_j.
\qquad (1\le i\le r).$$
Here $y_i=(y_{i1},...,y_{im})^t.$
The remaining left singular vectors $u_{r+1},\ldots,  u_m$ can be chosen to be any orthonormal collection  of
$m-r$ vectors in $\C^N$ that are perpendicular to $\text{ span } \{a_1,\ldots, a_m\} $.   One may then readily verify that (\ref{SVD}) holds.

 The  {\em Frobenius norm} of a matrix $X=[x_1,\ldots, x_m]\in \C^{N\times m}$
 is $\|X\|_F = tr(X^*X),$ where {\em tr } denote the  trace of a matrix.

Now, the following approximation theorem of Schmidt (cf. \cite{Sch07}) and later rediscovered by Eckart and Young (\cite{EY36}) shows that the 
 SVD can be used to find  the subspace of dimension $n$ that is closest to a given finite numbers of vectors.  

\begin{theorem} \label{EY}
Let $\{a_1, \dots, a_m\}$, be a set of vectors in $\C^N,$ such that $r =$ dim (span$\{a_1, \dots, a_m\})$, and 
suppose  $A =[a_1,\ldots, a_m]$, has SVD $A=U\Lambda^{1/2} Y^*$ with $0<n\le r$.  Then 
$A_n:=\sum_{j=1}^n \sqrt{\lambda_j} u_j y_j^*$   satisfies
$$
\|A-A_n\|_F=\min_{\text{rank }B\le n}\|A-B\|_F=\left(\sum_{j=n+1}^r \lambda_j\right)^{1/2}.
$$
If $\lambda_{n+1}\neq \lambda_n$, then $A_n$ is the unique such matrix of rank at most $n$. 
\end{theorem} 
Equivalently, 
\begin {theorem} \label{E-Y}
Let $\{a_1, \dots, a_m\}$, be a set of vectors in $\C^N$ such that $r=$ dim (span$\{a_1, \dots, a_m\}) $, and 
suppose  $A =[a_1,\ldots, a_m]$, has SVD $A=U\Lambda^{1/2} Y^*$ and that $0<n\le r$.  
 If  $W=\sspan \{u_1,\dots,u_n\}$, then
$$\{P_{W}a_1,\ldots, P_{W}a_m\} =\sum_{i=1}^n \sqrt{\lambda_i}u_iy_i^*   =A_n$$ and
\begin {equation}
\sum_{i=1}^m \|a_i-P_{W}a_i\|_{2}^2 \le \sum_{i=1}^m \|a_i-P_{M}a_i\|_{2}^2, \quad \forall \; M,  \dim M\le n,
\end {equation}
and the   space $W$ is unique if $\lambda_{n+1}\neq \lambda_n$.
\end {theorem}

\section*{Acknowledgements}
We thank Yves Meyer for the encouragement to write this paper, and the anonymous referees for helpful suggestions which improved the final version of this manuscript.


\providecommand{\bysame}{\leavevmode\hbox to3em{\hrulefill}\thinspace}
\providecommand{\MR}{\relax\ifhmode\unskip\space\fi MR }
\providecommand{\MRhref}[2]{%
  \href{http://www.ams.org/mathscinet-getitem?mr=#1}{#2}
}
\providecommand{\href}[2]{#2}

\end{document}